\RequirePackage{fix-cm}
\documentclass[smallextended]{svjour3}       
\usepackage{graphicx}
\usepackage[applemac]{inputenc}		
\usepackage{lmodern}						
\usepackage[english]{babel}					
\usepackage[weather]{ifsym}	
\usepackage{amsmath}  
\usepackage{amsfonts}  
\usepackage{amssymb}  

\usepackage{ulem}

\usepackage{latexsym}

 \journalname{Journal of Statistical Physics}

\numberwithin{equation}{section}

\def\R {\mathbb{R}}
\def\D {\mathbb{D}}
\def\T {\mathbb{T}}

\def\cA {\mathcal{A}}

\def\cC {\mathcal{C}}

\def\cD {\mathcal{D}}
\def\cK {\mathcal{K}}

\def\cH {\mathcal{H}}

\def\cG {\mathcal{G}}
\def\cL {\mathcal{L}}
\def\cF {\mathcal{F}}
\def\cM {\mathcal{M}}
\def\cN {\mathcal{N}}

\def\cJ {\mathcal{J}}

\def\cP {\mathcal{P}}

\def\cZ {\mathcal{Z}}

\def\cR {\mathcal{R}}

\newcommand{\bbB}{{\ensuremath{\mathbb B}} }

\newcommand{\bbD}{{\ensuremath{\mathbb D}} }
\newcommand{\bbE}{{\ensuremath{\mathbb E}} }

\newcommand{\bbP}{{\ensuremath{\mathbb P}} }

\newcommand{\bbR}{{\ensuremath{\mathbb R}} }

\newcommand{\bbT}{{\ensuremath{\mathbb T}} }

\newcommand{\gd}{\delta}
\newcommand{\gep}{\varepsilon}       
\newcommand{\gp}{\varphi}
\newcommand{\gr}{\rho}
\newcommand{\gz}{\zeta}

\newcommand{\gP}{\Phi}

\newcommand{\gl}{\lambda}
\newcommand{\gL}{\Lambda}
\newcommand{\gs}{\sigma}

\def\a {{\alpha}}
\def\b {{\beta}}
\def\g {{\gamma}}

\def\eps {{\varepsilon}}
\def\e {{\varepsilon}}

\def\indc {{\bf 1}}

\def\cum{{f^\eps_{n,[0,t]} (H ^{\otimes n } )}}
\def\lcum{{f_{n,[0,t]} (H^{\otimes n} )}}

\def\d {{\partial}}

\newcommand{\ba}{\begin{aligned}}
\newcommand{\ea}{\end{aligned}}

\newcommand{\be}{\begin{equation}}
\newcommand{\ee}{\end{equation}}

\newcommand{\mb}{  \Delta  \hskip-.2cm  \Delta }

\begin{document}

\title{Fluctuation theory in the Boltzmann--Grad limit}


\author{Thierry Bodineau \and Isabelle Gallagher \and Laure Saint-Raymond \and Sergio Simonella}
\institute{
T. Bodineau \at CMAP, CNRS, Ecole Polytechnique, I.P. Paris\\ 
\email{thierry.bodineau@polytechnique.edu}
\and 
I. Gallagher 
\at DMA, Ecole normale sup\'erieure \\
\email{isabelle.gallagher@ens.fr}
\and 
L. Saint-Raymond
\at UMPA, Ecole normale sup\'erieure de Lyon\\
\email{laure.saint-raymond@ens-lyon.fr}
\and
S. Simonella
\at UMPA, Ecole normale sup\'erieure de Lyon\\
\email{sergio.simonella@ens-lyon.fr}
}


\date{Received: date / Accepted: date}

\maketitle

\begin{abstract}
We develop a rigorous theory of hard-sphere dynamics in the kinetic regime, away from thermal equilibrium. 
In the low density limit, the empirical density obeys a law of large numbers
and the dynamics is governed by the Boltzmann equation.  Deviations from this behaviour
are described by   dynamical correlations, which can be fully
characterized for short times. This provides both a fluctuating Boltzmann equation 
and large deviation asymptotics.
\end{abstract}

 \section{Introduction}

 In this note we report some recent progress on the origins of the fluctuation theory 
 from the fundamental laws of motion. For states far from equilibrium, the macroscopic fluctuation theory has been  investigated intensively, but microscopic derivations are mainly focused on stochastic lattice gases
 (see e.g.~\cite{S2,BdSGJ-LL,derrida}). 
 We study here classical deterministic particles in a rarefied gas. In its rigorous version,
 the issue is then connected with the problem of the mathematical validity of the Boltzmann equation,
 in the limit introduced by Grad~\cite{Gr49}. This limit procedure states that in a Hamiltonian system of $N$ particles,
 strongly interacting at distance $\e$, the  particle density approximates the solution to the Boltzmann equation when $N \to \infty ,$ $\e \to 0 ,$ in such a way that the collision frequency (proportional to $N\e^{d-1}$ in dimension $d=2,3$) remains bounded; the volume density scales like $\e$, and both collisions and transport have a finite effect in the limit.
 
The Boltzmann gas is a simple case featuring a nonlinear dynamics, and a rich structure for the fluctuations.
At the macroscopic scale, typical particles behave as i.i.d. variables. Small fluctuations
admit already an interesting theory. In particular, they exhibit  spatial correlations, and noise originating from (deterministic) collisions  \cite{EC81,S83}. 
Moreover, rare fluctuations satisfy a large deviations principle. We refer to the companion work
by F. Bouchet \cite{bouchet}, where the large deviations theory for the Boltzmann equation has been discussed first. 
 
On the mathematical side, the only result we are aware of in this direction, is the convergence of the second
moment of the small fluctuations proved by H. Spohn in \cite{S81}. Similar results are available
for linear regimes close to equilibrium, both for short \cite{BLLS80} and for large times \cite{BGSR2}.
As suggested in \cite{S83}, the fluctuation theory should  not be merely a phenomenological theory, but a 
rigorous consequence of the laws of mechanics. Our aim is to support this assertion, providing a 
robust mathematical framework. 

We shall state several theorems (Theorems \ref{LanfordThm}, \ref{thmTCL}, \ref{thmLD} below)
describing the behaviour of the  empirical density 
\begin{equation}
\label{eq: empirical}
\pi^{\e}_t:=\frac{1}{\mu_\e}\sum_{i=1}^N \gd_{{\bf z}^\e_i(t)} \,,\ \ \ \ \ \mu_\e =\e^{-(d-1)}
\end{equation}
for a Newtonian evolution of $N$ particles with positions and velocities $${\bf z}^\e_i(t) = \left({\bf x}^\e_i(t), {\bf v}^\e_i(t)\right)\;,\ \ \ \ \  i=1,\dots,N\;.$$ 
We assume that the particles are approximately Poisson-distributed 
at time~$t=0$, with (random) total number of particles $\cN$ and
regular phase-space density~$f^0=f^0(x,v)$. Probability and expectation with
respect to this initial measure are denoted by $\bbP_\e$ and $\bbE_\e$.
Then, in the Boltzmann-Grad limit $\e \to 0$, $\bbE_\e\left(\cN\right)/\mu_\e \to 1$, the following properties hold.

\begin{enumerate}
\item {\it Law of large numbers:} 
\begin{equation} 
\label{eq:R1}
\pi^{\e}_t \to f_t  \,,\qquad  t \in [0,T^\star]
\end{equation}
 weakly (in probability) for some $T^\star>0$, where $f_t$ is the solution of
 \begin{equation} 
 \label{eq:Boltz}
\d_t f +v \cdot \nabla _x f = C(f,f)
\end{equation}
with initial datum $f^0$ and $C$ is Boltzmann's collision operator  \cite{La75}; in particular, the {\it chaos} property of the initial measure propagates in time (rescaled correlation functions converge to a tensor product).

\medskip
\item {\it Central limit theorem:} the fluctuation field  
\begin{equation}
\label{eq: fluctuation field}
\gz^\gep_t := \sqrt{\mu_\eps }
\left( \pi^\e_t -  \bbE_\eps\left(  \pi^\e_t  \right) \right)
\end{equation} 
describing the small deviations of the empirical density 
 from its average,
  converges in law  on $[0,T^\star]$ to the Gaussian process $\gz_t$ governed by the fluctuating Boltzmann equation
\begin{equation} \label{eq:R2}
d \gz_t   = \cL_t \,\gz_t\, dt + d\eta_t  \,,
\end{equation}
where $\cL_t$ is Boltzmann's operator linearized around $f_t$, and $d\eta_t$ is Gaussian noise
 (with covariance defined in \eqref{eq: covariance bruit} below),
as predicted in \cite{S81}.

\medskip
\item {\it Large deviations} are exponentially small in $\mu_\e$ and characterized, at least in a regime of strong regularity, by the same large deviation functional as heuristically derived in \cite{bouchet} (and previously obtained, rigorously, in \cite{Rez2} from a one-dimensional stochastic process). That is, the probability of observing a path $ \varphi_t=  \varphi(t,x,v)$ satisfies
\begin{equation} \label{eq:R3}
\bbP_{\eps}\left( \pi^\eps_t \approx  \varphi_t , \ t \in [0, T^\star] \right) \asymp
\exp\left(-\mu_\e\, \cF(T^\star,\varphi)\right) ,
\end{equation}
where $\cF$ is defined as the Legendre transform of a functional $\cJ = \cJ( T^\star,\varphi)$, solution of a Hamilton-Jacobi 
equation.
\end{enumerate}

The Boltzmann equation is naturally suited to a probabilistic interpretation, and its mathematical validity can be based on the construction of a stochastic particle system mimicking the microscopic collisions.
The basic example is the Kac model \cite{Kac}, from which the spatially homogeneous Boltzmann equation can indeed be recovered. Fluctuations in this type of processes can also be analysed (see e.g. \cite{KL,vK74,meleard}), including large deviations (\cite{Le75}) and spatially inhomogeneous variants (\cite{Rez,Rez2}).
Our results show that the analogy between the deterministic hard-sphere dynamics and the stochastic model 
goes far beyond the typical behavior as it remains valid for extremely rare events. 
The statistical behavior of the hard-sphere gas described above is in fact the same as the one derived in 
\cite{Rez,Rez2}.
For physical observables as the empirical measure, the deterministic dynamics and the stochastic approximation (as often used in simulations) cannot be distinguished,
even at the level of fluctuations.

Our main restriction is  the smallness of the time $T^\star$. This time (depending only on $f^0$) is actually  a fraction of the
time of validity of the Boltzmann equation obtained by Lanford in \cite{La75}. We will
further restrict to a gas of hard spheres, though we believe that the results could be proved for
smooth and compactly supported interactions, adopting known techniques~\cite{Ki75,GSRT,PSS17}.

The Hamilton-Jacobi equation determining $\cJ$ (Theorem \ref{thm: HJ} below) is our ultimate point of arrival
in the derivation from a microscopic mechanical model. A stationary solution of this
equation is given by the (dual) Boltzmann's $H$ functional, which describes large deviations of the
equilibrium state. Moreover,~$\cF$ has an invariance encoding the microscopic reversibility,
a symmetry inherited from the equality between the probability of a path and the probability of the time-reversed 
path.  This is an indication on the amount of recovered information which was “lost” in Lanford’s Theorem, proving the transition from a reversible to a dissipative model.

Our method is far from standard approaches to a large deviation problem. 
For stochastic dynamics, large deviations can be evaluated by modifying the underlying stochastic process in a time dependent way in order to produce an atypical trajectory. The optimal cost for inducing such a bias  on the stochastic dynamics is precisely the large deviation rate. 
For hard-sphere systems, there is no underlying stochastic dynamics 
as all the randomness 
lies in the initial data. It seems exceedingly
hard to figure out a way to bias the initial probability measure in order to produce a given admissible path $\varphi_t$. Indeed, 
the deterministic dynamics is responsible for an intricate relation between the path and the initial distribution of spheres.

We therefore turn back to the more modest problem of analysing the error in \eqref{eq:R1}\footnote{For previous quantitative investigations of the correlation error, we refer to \cite{GSRT,PS17}.}.
It is already evident in Lanford's proof, that the dynamical information
lives on precise little regions of the $j-$particle phase space, converging
to measure-zero sets as $\e \to 0$, for any finite $j$. In  little regions of the same size, correlations are generated by the
collision events, which break the propagation of chaos (see e.g.\,\cite{BGSRS18}). These correlation sets do not encode the most probable future dynamics and they can be neglected when proving \eqref{eq:R1}.
However we can extract much more information, by looking for mathematically
tractable quantities which are concentrated exactly on these sets, and retaining the 
information which is lost in \eqref{eq:R1}. 

A natural candidate is provided by {\it cumulants}, which can be obtained by the series expansion of
the generating function
\begin{equation}\label{def-Lambdaepst}
\gL^\eps_{t} (h) 
:= \frac{1}{\mu_\eps}\log \bbE_\eps \left( \exp \left(  \sum_{i =1}^{\cN} h \left( {\bf z}^\e_i(t)\right)  \right)   \right)
\end{equation}
where $h$ is any test function. The order $n$ in this expansion is given by a function
$f_n^{\e}= f_n^{\e}(t)$ on the $n-$particle phase space (formula \eqref{exponential moment} below) describing a cluster of particles mutually correlated by a chain of interactions.
It has been noted in \cite{EC81} that the hierarchy of cumulants
determines {\it all} the properties of the fluctuations in a gas,  and that their exact computation
furnishes a theory of fluctuations at the same time. In order to prove rigorous results and reach
large deviations, we will construct the limit of the exponential moment \eqref{def-Lambdaepst},
and link it to the function $\cJ$. 

The expansion of \eqref{def-Lambdaepst} leads to a combinatorial problem,
which can be dealt with by the cluster expansion method \cite{Ru69}. Indeed this method fits very well with the dynamics at low
density, when combined with geometrical estimates on hard-sphere trajectories.

We organize the paper as follows.
Section \ref{strategy} is a brief introduction to our strategy. Section~\ref{sec:RPT} presents the model and the fundamental result leading to~\eqref{eq:R1}, and explains the basic dynamical formula expressing the main quantities of interest in terms of the initial data. In Section \ref{sec: Dynamical correlations} we state our main results on the dynamical correlations and their limiting structure, and  derive the Hamilton-Jacobi equation. Finally, the last two sections are devoted to the fluctuating Boltzmann equation and the large deviations respectively. In this paper we shall only sketch the proof of our results, the complete version of which will be provided in a longer publication \cite{BGSRS}.

\section{Strategy} \label{strategy}

Lanford's method \cite{La75} is based on the BBGKY hierarchy governing the evolution of the family of (properly rescaled) correlation functions~$\left(F^{\eps}_n \right)_{n \geq 1}$. This hierarchy is completely equivalent to the Liouville equation describing interacting transport of $N$ hard spheres. In the Boltzmann-Grad limit, the probability densities concentrate on an infinite-dimensional space, and the BBGKY hierarchy  is convenient to capture the relevant information. One thus introduces the (rescaled)  correlation functions $ F_n^{\eps } (t,Z_n )$ such that
$$\bbE_\eps \Big( \sum_{\substack{i_1, \dots, i_n \\ i_j \neq i_k, j \neq k}} h_n \big( {\bf z}_{i_1}^{\e }(t), \dots ,  {\bf z}_{i_n}^{\e } (t)\big) \Big)
=  \mu_\eps^{n} \int_{\D^n} dZ_n \, F_n^{\eps } (t,Z_n )\, h_n( Z_n)   \,,
$$
for any test function $h_n$. The family $\left(F^{\eps}_n \right)_{n \geq 1}$ is suited to the description of typical events: in the limit $\e \to 0$, $F^{\eps}_n\to f^{\otimes n}$ so that everything is coded in $f$ (solution of \eqref{eq:Boltz}), no matter how large~$n$. 

We need to go beyond  the BBGKY hierarchy and turn to  a more powerful representation of the dynamics.
We shall replace the family $\left(F^{\eps}_n \right)_{n \geq 1}$ with an equivalent family of (rescaled) truncated correlation functions $\left(f^{\eps}_n \right)_{n \geq 1}$,  called  {\it cumulants}. Their role is to  grasp information on the dynamics on finer and finer scales. Loosely speaking, $f^{\eps}_n $ will collect events where~$n$ particles are ``completely connected'' by a chain of interactions. We shall say that the $n$ particles form a connected  {\it cluster}. Since a collision between two given particles is typically of order~$\mu_\eps^{-1}$ (the size of the ``collision tube'' spanned by one particle in time 1), a complete connection would account for events of probability  of order~$\mu_\eps^{-(n-1)}$. We therefore end up with a hierarchy of rare events, which we would like to control at arbitrary order. At variance with $\left(F^{\eps}_n \right)_{n \geq 1}$, even  {\it after} the limit $\mu_\eps \to \infty$ is taken, the cumulant $f^{\eps}_n$ cannot be trivially obtained from the cumulant $f^{\eps}_{n-1}$. Each step entails extra information, and events of increasing complexity, and decreasing probability.

Unfortunately, the equations for $\left(f^{\eps}_n \right)_{n \geq 1}$ are difficult to handle. But the moment-to-cumulant relation~$ \left(F^{\eps}_n\right)_{n \geq 1} \to  \left(f^{\eps}_n \right)_{n \geq 1}$ is a bijection and, in order to construct $ f^{\eps}_n(t)$, we can still resort to the same solution representation of~\cite{La75} for the correlation functions $\left(F^{\eps}_n (t)\right)_{n \geq 1}$.
 This formula is an expansion over  {\it collision trees}, meaning that it has a geometrical representation as a sum over binary tree graphs, with vertices accounting for collisions (see Section~\ref{sec:RPT}). Two particles 
are correlated if their generated trees are  {\it connected} by a ``recollision'', which is an event of weight $\mu_\eps^{-1}$
(see Section \ref{loops} for a precise notion of recollision).

In Proposition \ref{prop:APEC} we will state the main technical advance of this paper: the  cumulant (rescaled by the factor $\mu_\eps^{n-1}$) grows as $ n^{n-2}$ in $L^1$-norm. This estimate is intuitively simple. We have at disposal a geometric notion of correlation as a {\it link} between two collision trees.
Based on this notion, we can draw a random graph on $n$ vertices telling us which particles are correlated and which particles are not (each collision tree being one vertex of the graph). Since the cumulant $f^\eps_n$ corresponds to $n$ completely correlated particles, there will be at least $n-1$ edges, each one of small `volume' $\mu_\eps^{-1}$. Of course there  could  be more than $n-1$ connections (the random graph has cycles), but these are hopefully unlikely as they produce extra smallness in $\eps$. If we ignore all of them, we are left with minimally connected graphs, whose total number is $n^{n-2}$ by Cayley's formula. 

The limiting equations for the family $\left(f^{\eps}_n \right)_{n \geq 1}$ form a {\it Boltzmann cumulant hierarchy}, displaying a remarkable structure \cite{EC81}. The first equation ($n=1$) is just the Boltzmann equation. The second equation ($n=2$) is driven by a linearized Boltzmann operator $\cL_t$, plus a singular ``recollision operator'', acting on $f_1$ only, generating the ``connection'' (correlation) between two particles and suited to be interpreted as noise source \cite{S81,cohen}. The higher order equations ($n>2$) have an increasingly complex structure, combining the action of the two operators (of standard linearized type, and of connecting type) on $n$ different particles, in all possible ways. But  the good $n$-dependence of the uniform bounds allows to sum up the cumulants into an analytic series. This finally translates the cumulant hierarchy into the Hamilton-Jacobi equation, which stands as a compact, nonlinear representation of the correlation dynamics.


 \section{Collision trees} \label{sec:RPT}

In this section we introduce the geometrical representation of the hard-sphere dynamics with random initial data, which will be our basic tool.

  \subsection{Hard-sphere model}

The microscopic model consists of $N$ identical hard spheres of unit mass and of diameter~$\eps$. \begin{figure}[h] 
\centering
\includegraphics[width=2.5in]{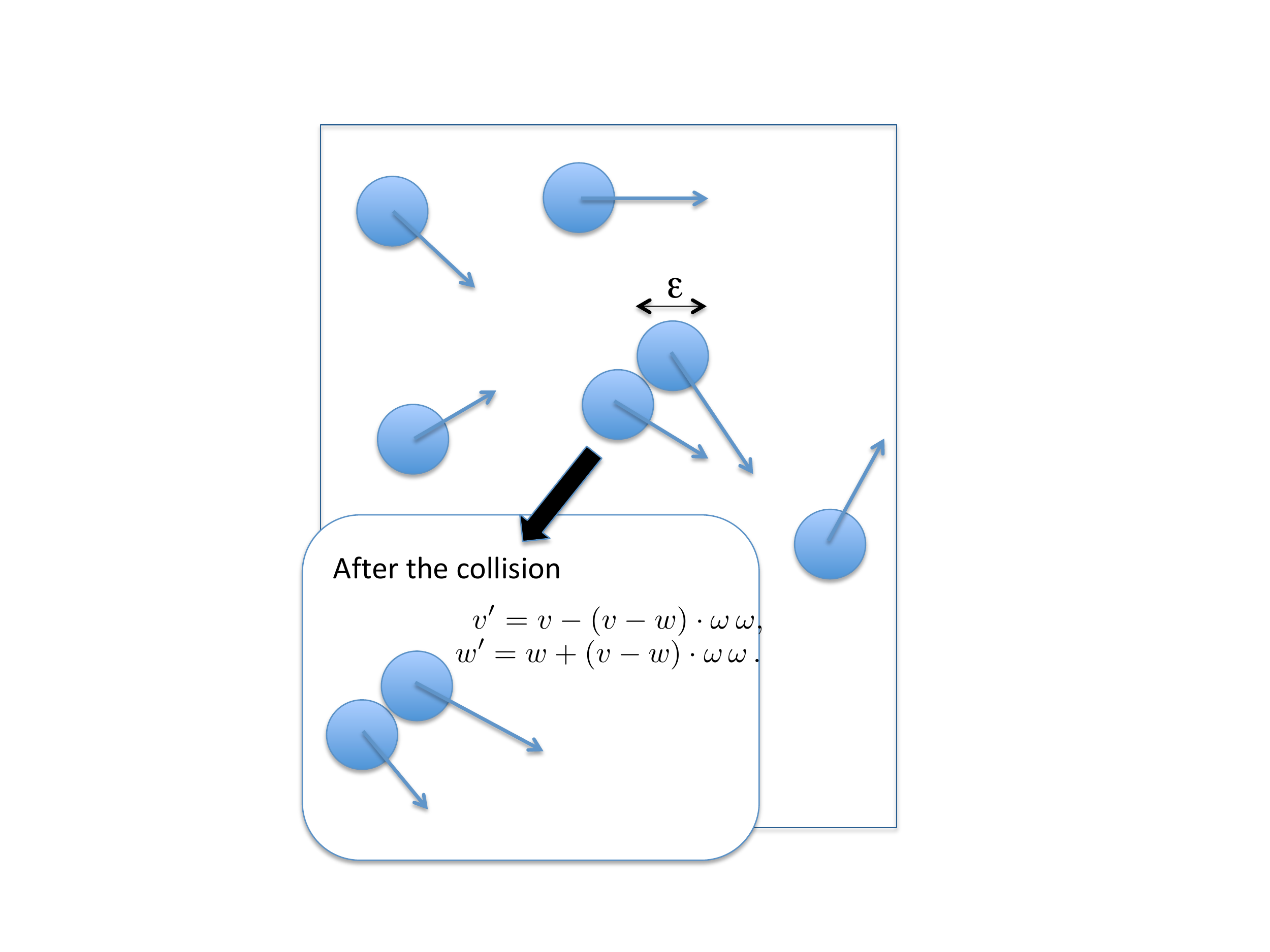} 
\caption{\small 
Transport and collisions in a hard-sphere gas
}
\end{figure}
Their motion is governed by a system of ordinary differential equations, which are set in~~$\D^N:=( \T ^d\times\R^d)^{N }$ where~$\mathbb T^d$ is  the unit $d$-dimensional periodic box:
$$
{d{\bf x}^{\e}_i\over dt} =  {\bf v}^{\e}_i\,,\quad {d{\bf v}^{\e}_i\over dt} =0 \quad \hbox{ as long as \ } |{\bf x}^{\e}_i(t)-{\bf x}^{\e}_j(t)|>\eps  
\quad \hbox{for \ } 1 \leq i \neq j \leq N
\, ,
$$
with specular reflection at collisions: 
\begin{equation}
\label{defZ'nij}
\begin{aligned}
\left. \begin{aligned}
 \left({\bf v}^{\e}_i\right)'& := {\bf v}^{\e}_i - \frac1{\eps^2} ({\bf v}^{\e}_i-{\bf v}^{\e}_j)\cdot ({\bf x}^{\e}_i-{\bf x}^{\e}_j) \, ({\bf x}^{\e}_i-{\bf x}^{\e}_j)   \\
\left({\bf v}^{\e}_j\right)'& := {\bf v}^{\e}_j + \frac1{\eps^2} ({\bf v}^{\e}_i-{\bf v}^{\e}_j)\cdot ({\bf x}^{\e}_i-{\bf x}^{\e}_j) \, ({\bf x}^{\e}_i-{\bf x}^{\e}_j)  
\end{aligned}\right\} 
\quad  \hbox{ if } |{\bf x}^{\e}_i(t)-{\bf x}^{\e}_j(t)|=\eps\,.
\end{aligned}
\end{equation}
The sign of the scalar product $({\bf v}^{\e}_i-{\bf v}^{\e}_j)\cdot ({\bf x}^{\e}_i-{\bf x}^{\e}_j)$ identifies post-collisional (+) and pre-collisional (-) configurations.
This flow does not cover all possible situations, as multiple collisions are excluded.
But one can show (see \cite{Ale75}) that, for almost every initial configuration $({\bf x}^{\e 0}_i, {\bf v}^{\e 0}_i)_{1\leq i \leq N}$, there are neither multiple collisions,
nor accumulations of collision times, so that the dynamics is globally well defined. 

Below, we shall denote collections of positions and velocities respectively by~$X_N := (x_1,\dots,x_N) \in \T^{dN}$ and~$V_N := (v_1,\dots,v_N) \in \R^{d N}$, and we set~$Z_N:= (X_N,V_N) \in ( \T ^d\times\R^d)^{N }$, $Z_N =(z_1,\dots,z_N)$.

 Let $f^0$ be a probability density on $\D $ with Gaussian decay in velocity
\begin{equation}
\label{beta0}
 |f^0(x,v)|  + | \nabla_x f^0(x,v)|\ \leq C_0 \; \exp  \left(- \frac{\beta_0} 2 |v|^2 \right) , 
 \end{equation}
where $C_0, \beta_0 > 0$.
Because of the condition of hard-sphere exclusion, the positions of the particles cannot be independent of each other.
To better focus on the dynamical issue, we shall choose, as initial measure, the $N$-particle distribution with minimal correlations.
%
%
In particular, to avoid spurious correlations due to a given total number of particles, we shall consider a grand canonical state. 
The initial probability density of finding $N$ particles in $Z_N$ is given by
\begin{equation}
\label{eq: initial measure}
\frac{1}{N!}W^{\eps 0}_{N}(Z_N) 
:= \frac{1}{\cZ^ \eps} \,\frac{\mu_\eps^N}{N!} \, \indc_{
{\mathcal D}^{\eps}_{N}}
 \, \prod_{i=1}^N f^0 (z_i)  \end{equation} 
where the domain encodes the exclusion:
$${\mathcal D}^{\eps}_{N}
:= \big \{ Z_N \in  {\mathbb D}^N \, \big | \, \quad  \forall i \neq j, \quad  \, |x_i - x_j| >\eps \big\}\,,
$$
and the normalization constant $\cZ^\eps$ is given by 
$$
\cZ^\eps :=  1 + \sum_{N\geq 1}\frac{\mu_\eps^N}{N!}  
\int_{\D^N} dZ_N\, \indc_{
{\mathcal D}^{\eps}_{N}} \prod_{i=1}^N f^0 (z_i) \,  .
$$
With this definition, if $$\mu_\eps \eps^{d-1} = 1\;,$$ then
the average number of particles satisfies
$$
\lim_{\gep \to 0} {\mathbb E}_\eps  \left(  \cN  \right)  \eps^{  d-1 } =  1
$$
(Boltzmann-Grad scaling).

The rescaled {\it $n$-particle correlation function} is defined by
\begin{equation}\label{defFeps0n}
\begin{aligned}
F_n^{\eps 0} (Z_n)
:= \mu_\eps^{-n} \,
\sum_{p=0}^{\infty} \,\frac{1}{p!}\, \int dz_{n+1}\dots dz_{n+p} \,
W_{n+p}^{\eps 0} (Z_{n+p})\;.
\end{aligned}
\end{equation}
For any symmetric test function $h_n : \D^n\rightarrow \bbR$, one can check that
\begin{equation}
\label{eq: marginal time 0}
\begin{aligned}
\bbE_\eps \Big( \sum_{\substack{i_1, \dots, i_n \\ i_j \neq i_k, j \neq k}} h_n \big( {\bf z}_{i_1}^{\e 0}, \dots ,  {\bf z}_{i_n}^{\e 0} \big) \Big)
=  \mu_\eps^{n} \int_{\D^n} dZ_n \, F_n^{\eps 0} (Z_n )\, h_n( Z_n)   \,.
\end{aligned}
\end{equation}
Moreover one can prove that, in the Boltzmann-Grad limit,
$$
\forall n \geq 1 \, , \quad  F_n^{\eps 0} (Z_n) \longrightarrow \prod_{i=1}^n f^0 (z_i) \hbox{ as } \eps \to 0$$
on the set $\{ x_i \neq x_j\;,\ \forall i \neq j\}$.
That is, at leading order, the initial distribution is  chaotic.

Starting from the dynamical equations (\ref{hardspheres}) we get that, for each fixed~$N$, the probability density at time $t>0$ is determined by the Liouville equation
\begin{equation}
\label{Liouville}
	\d_t W^{\eps}_N +V_N \cdot \nabla_{X_N} W^{\eps}_N =0  \,\,\,\,\,\,\,\,\, \hbox{on } \,\,\,{\mathcal D}^{\eps}_{N}\, ,
\end{equation}
	with specular reflection (\ref{defZ'nij}) on the boundary $|x_i - x_j|= \eps$.

  By integration of the Liouville equation for fixed $\eps$, we get that 
 the one-particle correlation function $F^\eps_1$   satisfies an equation
\begin{equation}
\label{first-correlation}
 \partial_t F^\eps_1 + v \cdot \nabla_{x} F^\eps_1 = C_{1,2}^{\eps} F^\eps_{2}  
 \end{equation}
where the collision operator comes from the boundary terms in Green's formula (using the reflection condition to rewrite the gain part in terms of pre-collisional velocities):
 $$
 \begin{aligned}
(C_{1,2} ^\eps F^\eps_2 )(x,v)
&:=  \int F^\eps_{2} (x,v', x+\e \omega,w') \big( (w- v)\cdot \omega \big)_+ \, d \omega dw\\
&\quad -   \int  F^\eps_{2} (x,v, x+\e \omega,w) \big( (
w - v
)\cdot \omega \big)_- \, d \omega dw\, ,
 \end{aligned}
$$
with
$$v' = v - (v-w) \cdot \omega\,  \omega, \quad w' = w +(v-w) \cdot \omega\,  \omega \,.$$
As in~(\ref{eq: marginal time 0}), $F_1^{\eps}(t)$ describes the average behavior of (identical) particles at time $t$:
$$
\bbE_\eps \left( \frac{1}{\mu_\eps} \sum_{i=1}^{\cN} h \left( {\bf z}^\e_i(t)\right) \right)=\int F^{\eps}_1(t,z)\, h(z)\, dz\,,
$$
for any test  function $h : \D\rightarrow \bbR$.  Similarly for any test function $h_2 : \D^2\rightarrow \bbR$,
the two-particle correlation function satisfies 
$$
\bbE_\eps \left( \frac{1}{\mu_\eps^2} \sum_{i\neq j  }  
h_2\left( {\bf z}_i(t),{\bf z}_j(t)\right) \right)=\int F^{\eps}_2(t,Z_2)\, h_2(Z_2)\, dZ_2\, .
$$

  \subsection{Law of large numbers}

 The issue with equation (\ref{first-correlation}) is that it is not closed:  it involves~$F^\eps_2$.
 At the level of (\ref{first-correlation}), Boltzmann's main assumption would correspond to the replacement
$$
\begin{aligned}
F^\eps_2 (t,x_1,v_1, x_2,v_2) \sim F^\eps_1 (t,x_1,v_1)F^\eps_1 (t,x_2,v_2)\, , 
\quad \mbox{as} \, \, \eps\to 0\,, \\ \mbox{when} \quad |x_1-x_2|=\e\;, \quad (x_1-x_2)\cdot(v_1-v_2)<0\, .
\end{aligned}$$
In other words, particles are assumed to be statistically independent, at least in pre-collisional configurations.
This very {\it strong chaos property} (which we assumed at time 0) is supposed to be valid for all times.

In the Boltzmann-Grad limit $\eps \to 0$, we then expect  $F^\eps_1$ to be well approximated by the solution to the {\it Boltzmann equation}
\begin{equation}
\label{boltzmann-eq}
 \d_t f +v \cdot \nabla _x f = C(f,f)
 \end{equation}
with 
$$ 
\begin{aligned}
 & C(f,f) (t,x,v)
 \\
 & \,  \,  \, := \int_{\R^d}\int_{{\mathbb S}^{d-1}} \Big(  f(t,x,w') f(t,x,v') - f(t,x,w) f(t,x,v)\Big) 
  ((v-w)\cdot \omega)_+ \, d\omega \,dw \,.
  \end{aligned}
$$

The claim by Boltzmann  that the  particle density is well approximated by Equation~(\ref{boltzmann-eq}) has been made rigorous by Lanford for short times.
%

\begin{theorem}[{\bf Lanford}, \cite{La75}]
\label{LanfordThm}
Consider a  gas  of hard spheres initially distributed according to  {\rm(\ref{eq: initial measure})}.
      Then, in the Boltzmann-Grad limit $\mu_\eps \to \infty$ with $\mu_\eps \eps ^{d-1}=1$,
the 1-particle distribution $F_1^\eps$   converges, uniformly on compact sets, towards the solution $f$ of the 
Boltzmann equation~{\rm(\ref{boltzmann-eq})}
on a short time interval  $[0, T^\star]$ 
(where~$T^\star$ depends on  the initial distribution $f^0$ through $C_0,\beta_0$ in {\rm(\ref{beta0}))}.

Furthermore for each~$n$, the  $n$-particle correlation function $F^\eps_n(t)$ converges almost everywhere to~$f^{\otimes n}(t)$ on the same time interval. 
\end{theorem}

The propagation of chaos obtained in Lanford's theorem implies in particular  that  the empirical measure  $\pi_t^\eps$, defined by (\ref{eq: empirical}), concentrates on the solution to the Boltzmann equation.
Indeed computing the variance we get, for any test function $h$, that
\begin{equation}\begin{aligned}
\label{Eepsh2} 
& \bbE_\eps \Big( \big(\pi^\eps_t(h)- \int F^\eps_1(t,z)\, h(z)\, dz \big)^2 \Big)\\
&   
=   \bbE_\eps \Big( \frac{1}{\mu_\eps^2}  \sum_{i=1}^{\cN} h^2 \big( {\bf z}^{\eps}_i(t)\big) 
+ \frac{1}{\mu_\eps^2}  \sum_{i \not = j}  h \big( {\bf z}^{\eps}_i(t)\big) h \big( {\bf z}^{\eps}_j(t)\big)
\Big) -   \Big(\int F^\eps_1(t,z)\, h (z)\, dz \Big)^2 \\
&  
=  \frac{1}{\mu_\eps}  \int F^\eps_1 \,h^2\, d z_1  
+ \int F^\eps_2\, h^{\otimes 2}\, d Z_2 -   \Big(\int F^\eps_1\, h\, dz \Big)^2
\end{aligned}
\end{equation}
which converges to 0 as $\eps \to 0$ since $F^\eps_2$ converges to $f^{\otimes 2}$ and $F^\eps_1$ to $f$.
This computation can be interpreted as a law of large numbers.

%
%
Theorem   \ref{LanfordThm} entails a drastic loss of information, which (as becomes clear from the proof) is retained in particular “recollision sets” of measure zero. Some of the microscopic time-reversible structure can be recovered by looking at correlations on finer scales. 
This is the role played by the rescaled dynamical {\it cumulants},
defined by
\begin{equation}
\label{cumulant-expansion-dyn}
f^\eps_n(t, Z_n):=  \mu_\eps^{n-1} \sum_{s= 1}^n  \sum _{\sigma \in \cP^s_n } (-1) ^{s-1}  (s-1) ! \, \prod_{i=1} ^s F^\eps_{|\gs_i|}(t,Z_{\gs_i})
 \,.
 \end{equation}
Here we denoted by $\cP^s_n$ the set of partitions of $\{1,\dots, n \}$ in $s$ parts, $\cP^s_n\ni\sigma = \left(\sigma_1,\dots,\sigma_s\right)$,  by $|\sigma_i|$ the cardinality of the set $\sigma_i$ and by $Z_{\sigma_i} = (z_j)_{j \in \sigma_i}$.
This formula is cooked up to extract the effect of recollisions and therefore to obtain the detailed correlation structure at arbitrarily small scales.
Note that, for fixed $\eps>0$, $\left(F^{\eps}_n \right)_{n \geq 1}$ and $\left(f^{\eps}_n \right)_{n \geq 1}$ provide the same amount of information, as shown by the inversion formula~:
\begin{equation}
\label{inversion}
F^{\e}_n(t,Z_n) = \sum_{s = 1}^{n} \sum_{\gs \in \cP^s_n}   \mu_\eps^{ - (n-s)} \prod_{i=1}^s f^{\eps}_{|\gs_i|}(t,Z_{\gs_i})\,.
\end{equation}
%

Before passing to the investigation to \eqref{cumulant-expansion-dyn}, we need to recall the main features of the proof of Theorem \ref{LanfordThm} (see \cite{CIP,S2,GSRT} for more details).

\subsection{Hierarchy and pseudotrajectories}
    
  The starting point is the equation  (\ref{first-correlation}) for the 1-particle correlation function~$F^\eps_1$.
    In order to get a closed system, we write similar equations for all correlation functions $F_n^\eps$
\begin{equation}
\label{BBGKYGC}
 \partial_t F^\eps_n + V_n \cdot \nabla_{X_n} F^\eps_n = C_{n,n+1}^\eps F^\eps_{n+1} \quad \mbox{on} \quad {\mathcal D}^\eps_{n  }\;,
\end{equation}
 with specular boundary reflection as in \eqref{Liouville} \cite{Ce72}. As~$C^\eps_{1, 2} $ above, $C^\eps_{n, n+1}$ describes
 collisions between one ``fresh'' particle (labelled $n+1$) 
 and one given particle~$i\in \{1,\dots, n\}$.

We denote by~$S^\eps_n$\label{Sn-def}  the group associated with free transport in $\cD^\eps_n$ (with specular reflection at collisions). 
Iterating Duhamel's formula, we can express the solution as a sum of operators acting on the initial data~:
\begin{equation} \label{eq:seriesexp}
F^\eps _n  (t) =\sum_{m\geq0}    Q^\eps_{n,n+m}(t) F_{n+m}^{\eps 0} \, ,
\end{equation}
where we have defined for $t>0$
$$
\begin{aligned}
Q^\eps_{n,n+m}(t) F_{n+m}^{\eps 0 }  := \int_0^t \int_0^{t_{1}}\dots  \int_0^{t_{m-1}}  S^\eps_n(t-t_{ 1}) C^\eps_{n,n+1}  S^\eps_{n+1}(t_{1}-t_{2}) C^\eps_{n+1,n+2}   \\
\dots  S^\eps_{n+m}(t_{m})    F_{n+m}^{\eps 0} \: dt_{m} \dots d t_{1} 
\end{aligned}$$
and~$Q^\eps_{n,n}(t)F^{\eps0}_{n} := S^\eps _n(t)F_{n}^0$,  $Q^\eps_{n,n+m}(0)F^{\eps 0} _{n+m} := \gd_{m,0}F^{\eps0}_{n+m}$.

%

In the following, we shall label $1*, \dots ,n*$ the $n$ particles with configuration $Z_n$ at time~$t$, and $1, \dots, m$ the $m$ ``fresh'' particles which are added by the collision operators. The configuration of the particle labeled~$i*$ will be denoted indifferently~$z_i^*=(x_i^*,v_i^*)$ or~$z_{i*}=(x_{i*},v_{i*})$.

\subsubsection{The tree structure}  
  
Each term of the series expansion \eqref{eq:seriesexp} (after inserting the explicit definition of the collision operators) can be represented by a collision tree $a = (a_i) _{i = 1, \dots, m}$, 
which records the combinatorics of collisions~:  the colliding particles at time $t_i$ are $i$ and $a_i \in \{1^*, \dots ,n^*\} \cup \{1, \dots, i-1\}$. We define the set~$\cA_{n,   m}$ of all possible such trees. Note that~$|\cA_{n,  m}| = n(n+1) \dots (n+m-1)$. Note also that, graphically, $a \in \cA_{n,   m}$ is represented by $n$ binary tree graphs (below, we will call collision tree both $a \in \cA_{n,   m}$ and each of its $n$ components).


For all collision trees $a \in \cA_{n,   k}$ and  all parameters  $(t_i, \omega_i, v_i)_{i=1,\cdots,m}$ with $t_1>t_2> \dots> t_m$, one constructs  {\it pseudo-trajectories} on $[0,t]$ $$\Psi^{\eps}_{n,m} = \Psi^{\eps}_{n,m} \Big(Z_n^*, (a_i, t_i, \omega_i, v_i)_{i=1,\dots,m}\Big)$$
iteratively on $i=1,2,\dots,m$
(denoting by $Z_{n,m}(\tau ) =\big (Z_{n}^*(\tau ) ,Z_m(\tau)\big)$ the coordinates of particles at time~$\tau \leq t_m$):
\begin{itemize}
\item starting from $(Z_n^*)$ at time $t =: t_0$,
\item transporting all existing particles backward on $[t_{i}, t_{i-1}]$  (on ${\mathcal D}^\eps_{n +i-1}$ with specular reflection at collisions),
\item adding a new particle $i$  at time $t_i$, with position $x _{a_i} (t_i) +\eps \omega_i$ and velocity $v_i$,
\item and applying the scattering rule (\ref{defZ'nij}) if $\big( x _{a_i} (t_i), v_{a_i}(t_i^+), x_{a_i} (t_i)+\eps \omega_i, v_i\big) $ is a post-collisional configuration.
\end{itemize}
We discard non admissible parameters for which this procedure is ill-defined; in particular we exclude values of $\omega_i$ corresponding to an overlap of particles (two spheres at distance strictly smaller than $\e$).
In the following we denote by~$\cG^\eps_m(a, Z_n^* )$ the set of admissible parameters.

The following picture is an example of such flow (for  $n=1, m=4$).
\begin{figure}[h] 
\centering
\includegraphics[width=3in]{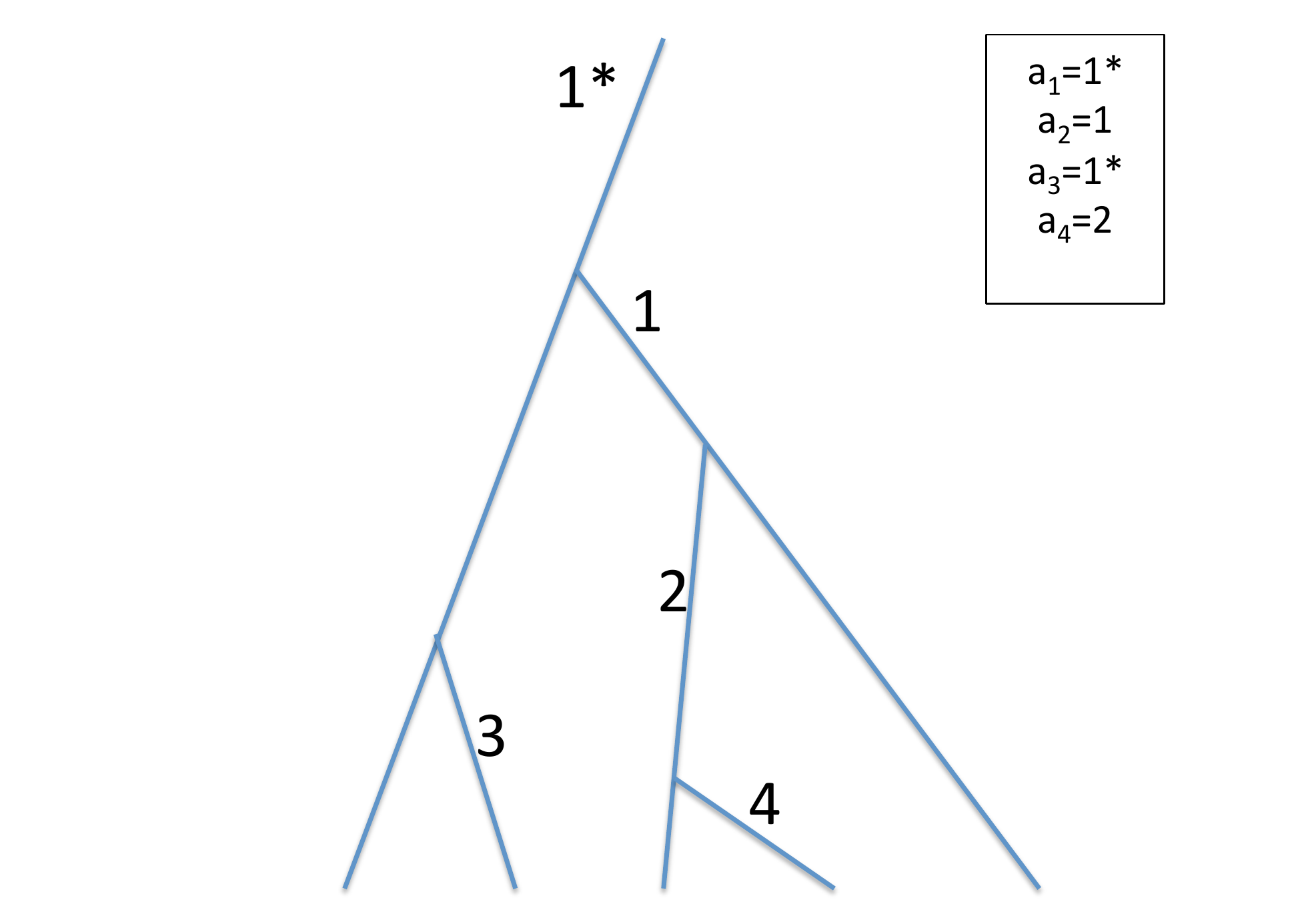} 
\caption{\small The tree structure of collisions.}
\end{figure}

With these notations, one gets the following geometric representation of 
the correlation function $F^\eps_n$ : 
$$
\begin{aligned}
 F^\eps _n (t, Z_n^*) =
 \sum_{m \geq 0} \sum_{a \in \cA_{n,   m} }\int_{\cG_{  m}^{\e}(a, Z_n^* )}    dT_m  d\Omega_{    m}  dV_{   m}\\
\times \left(\prod_{i=  1}^{  m}  \big( v_i -v_{  a_i} (t_i^+)\big) \cdot \omega_i \right)  
F_{n+m}^{\eps 0} \big (\Psi^{\eps0}_{n,m}\big)\;,
\end{aligned}
$$
where $(T_m, \Omega_{  m}, V_{  m}) := (t_i, \omega_i, v_i)_{  1\leq i\leq   m}$, and $\Psi^{\eps 0}_{n,m}$ is the $(n+m)$-particle configuration of the pseudo-trajectory at time zero.
Or, in short, 
\begin{equation}
\label{representation}
\begin{aligned}
 F_n^\eps(t,Z_n^*)&=  \int \mu (d  \Psi_{n}^\eps) \,  \cC \big(   \Psi_{n}^\eps \big) 
  \indc_{\cG^\e} (  \Psi_{n}^\eps \big)   F^{\eps 0} \big (   \Psi_{n}^{\eps 0}\big) \\
\hbox{ with } \quad\mu (d   \Psi_{n}^\eps ) &:=\sum_m\sum_{a \in \cA_{n,m} } dT_{m}  d\Omega_{m}  dV_{m},\quad 
\cC  (   \Psi_{n}^\eps  ) := \prod_{i=1}^{m}    \big( v_i -v_{a_i} (t_i^+)\big) \cdot \omega_i  
\end{aligned}
\end{equation}
$\indc_{\cG^{\e}} ( \Psi^\eps_n \big):=\indc_{\cG_{ m}^{\e}(a, Z^*_n)}$, and $F^{\eps 0} \big (\Psi_n^{\eps0}\big)$
the initial correlation function evaluated on the configuration at time 0 of 
the pseudo-trajectory (including $n+m$ particles). From now on, we will indicate by $\Psi^\e_n$ a generic pseudo-trajectory with~$n$ particles at time $t=t_0$.

 \subsubsection{A short time estimate}
 
 Each elementary integral corresponding to a collision tree with $m$ branching points involves a simplex in time ($t_1> t_2> \dots > t_m$). Thus, if we replace, for simplicity, the cross-section factors $\cC  (   \Psi_{1}^\eps  )$ by a bounded function (cutting off high energies), we immediately get that the integrals for $n=1$ are bounded, for each fixed tree~$a \in \cA_{1,m}$, by
$$\left|  \int   dT_m  d\Omega_{ m}  dV_{m}  \,  \cC( \Psi_{1}^\eps ) \, \indc_{\cG^\e} (  \Psi_{1}^\eps \big) \,
F^{\eps 0} \big ( \Psi_{1}^{\eps 0} \big)\right| \leq {(C'_0 t)^m \over m!}\,, $$
where $C'_0>0$ depends only on $C_0,\b_0$ of \eqref{beta0}. 
Since $|\cA_{1,m}| = m ! $, the series expansion is therefore absolutely convergent for short times, uniformly in $\e$. 
A similar estimate holds for $n > 1$. Moreover in    ce of the true factors~$\cC  (   \Psi_{n}^\eps  )$, the result remains valid (with a slightly different value of the convergence radius), though the proof requires some extra care \cite{Ki75}.

Hence it is enough to study the convergence  of each elementary term  in the Boltzmann-Grad limit $\eps \to 0$.

\subsubsection{Removing recollisions}\label{loops}

When the size $\eps$  of the particles goes to 0, we expect the pseudo-trajectory~$\Psi_{1}^\eps$ to converge to a limiting $\Psi_{1}$, defined iteratively on $i=1,2,\cdots,m$:

\begin{itemize}
\item starting from $z_1^*$ at time $t = t_0$,
\item transporting all existing particles backward on $[t_{i}, t_{i-1}]$ (by free transport),
\item adding a new particle $i$  at time $t_i$, exactly at position $x_{a_i} (t_i) $ and with velocity $v_i$,
\item  and applying the scattering rule (\ref{defZ'nij}) if $\big( v_i -v_{a_i} (t_i^+)\big) \cdot \omega_i  > 0$
(post-collisional configuration).
\end{itemize}
%

The main obstacle to the convergence $\Psi_{1}^\eps \to \Psi_{1}^\eps$ are the so-called recollisions. In the language of pseudo-trajectories, a {\it recollision} is a collision between pre-existing particles. Namely a collision which does not correspond to the addition of a fresh particle. It is easy to realize that, in the absence of recollisions,~$\Psi_{1}^\eps$ and $ \Psi_{1}$  differ only by small shifts in the positions.

 \begin{figure}[h] 
\centering
\includegraphics[width=3in]{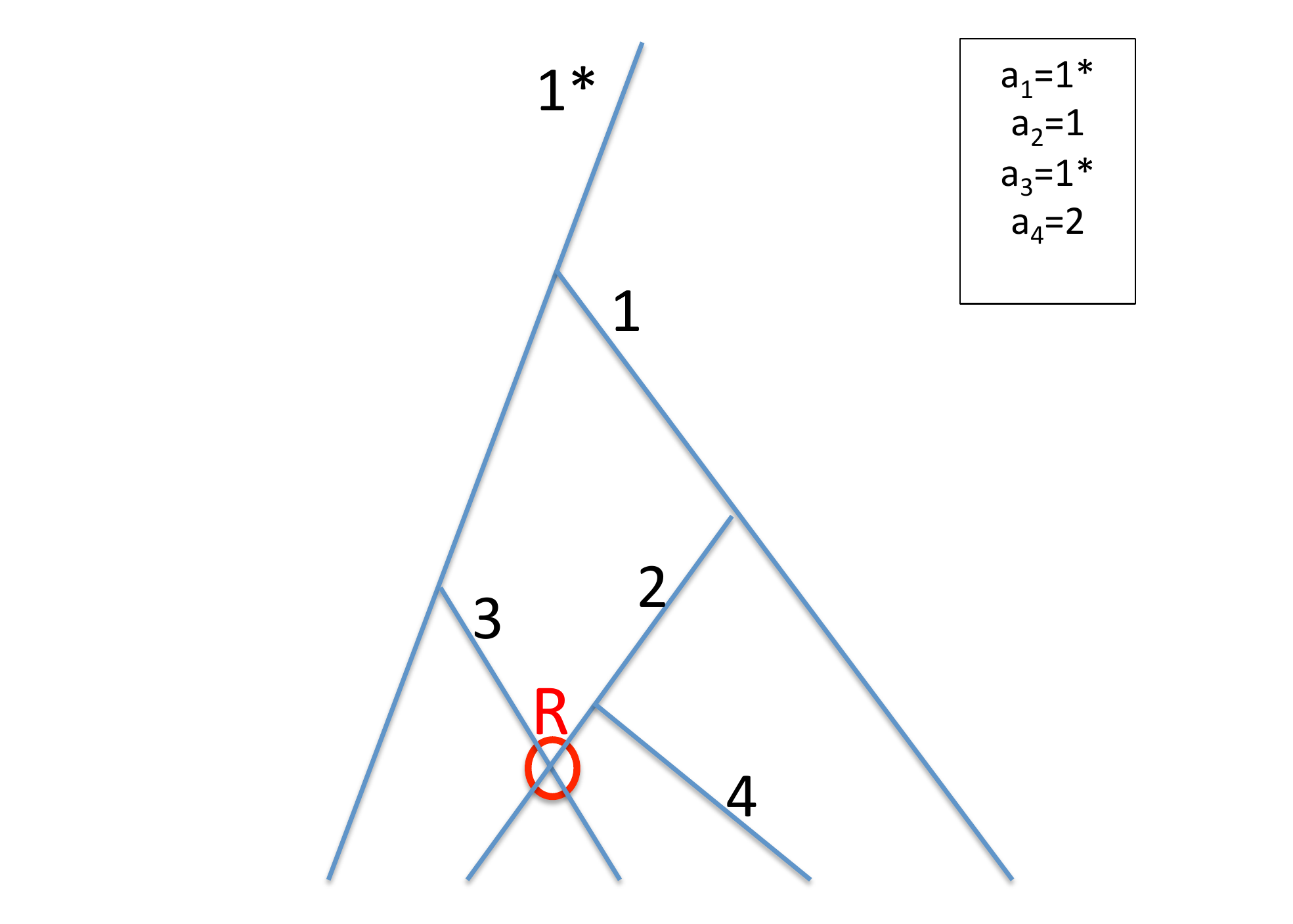} \label{recoll-fig}
\caption{\small 
An example of recollision. 
}
\end{figure}

A careful geometric analysis of recollisions shows that they can happen only for a small set of parameters, which is negligible in the limit~$\varepsilon \to 0$.  Roughly, if particles $p$ and $q$ are at positions $x_p, x_q$ with $x_p \neq x_q$ at  time $\tau > 0$, then 
a recollision between these particles implies  that there is a 
 time $t_{\rm rec} < \tau$ such that 
$x_p    -   x_q  - (v_p -v_q)  (\tau - t_{\rm rec}) = O(\eps)$. As a consequence, $v_p - v_q$ is constrained to be in a small cone of opening $\eps$, and the integration parameters in (\ref{representation}) lie in a small set.
Thanks to the uniform bounds, one concludes that pseudo-trajectories involving recollisions give an overall vanishing contribution to $F^\eps_1$.


%

A similar analysis can be performed to study higher order correlation functions. However, in this case the convergence is slightly more subtle since not all parameters are integrated. Notice that, for the $n$-particle correlation function, the convergence will fail on some sets of 
parameters of volume $\mu_\eps^{-1}$, which correspond to particles of {\it different} trees colliding in the backward dynamics (see e.g.\,Figure \ref{fig:rectreesf} below). These ``external recollisions" are apparently innocent, as they correspond again to small volume sets which do not contribute to the limit. On the other hand, it is the little failure of convergence of the $F_n^\eps$ which prevents  Lanford's theorem from being ``reversible", i.e. from being applicable to the state at time $t>0$ with reversed velocities \cite{BLLS80,BGSRS18}. This suggests that the relevant information to go backwards is hidden in singular directions where different trees merge. 
In Section \ref{sec: Dynamical correlations}, we will show that the dynamical cumulants $f_n^\eps$ defined by (\ref{cumulant-expansion-dyn}) ``live" in these singular directions, thus allowing to investigate the $n$-particle correlations in much higher detail.

 \subsubsection{Averaging over trajectories}

We conclude this section with a generalization of the previous discussion, which will be important in the sequel. 

So far we discussed correlations in phase space, at a given time $t$. But clearly, spatio-temporal correlations are also of interest. We therefore need to study trajectories of particles, and not only their distribution at a given time. 
Pseudo-trajectories provide a geometric representation of the iterated Duhamel series, but they are not physical trajectories of the particle system. Nevertheless,  the probability of trajectories of $n$ particles can be represented as above, by {\it conditioning} the Duhamel series.


\begin{proposition}[\cite{BGSRS}] \label{prop:APEC}
Let $H_n$ be a bounded measurable function on the Skorokhod space of trajectories over
$\D^n$ in $[0,t]$. Define
$$\begin{aligned}
 F^\eps _{n, [0,t]}  ( H_n) :=
 \int  dZ_n^* \int \mu (d  \Psi_{n}^\eps) \,  \cC \big(   \Psi_{n}^\eps \big) 
  \indc_{\cG^\e} (  \Psi_{n}^\eps \big)  
  H_n \big(Z_n^*([0,t]) \big) 
F^{\eps 0} \big (   \Psi_{n}^{\eps 0}\big)\,,
\end{aligned}
$$
where $Z_n^*( [0,t] )$ are the trajectories of the $n$ $*$-tagged particles in the pseudo-trajectory~$\Psi_{n}^\eps$. Then,
$$\begin{aligned}
\bbE_\eps \Big( \sum_{\substack{i_1, \dots, i_n \\ i_j \neq i_k, j \neq k}} H_n \big( {\mathbf z}^\eps_{i_1}([0,t] ), \dots ,  {\mathbf z}^\eps_{i_n}([0,t] ) \big) \Big)
=
\mu_\eps^{n} F^{\eps}_{n, [0,t]} (H_n)  \,,
\end{aligned}
$$
where ${\mathbf z}^\eps_{i_1}([0,t] ), \dots ,  {\mathbf z}^\eps_{i_n}([0,t] )$ is  
the sample path of $n$ hard spheres labeled $i_1, \dots, i_n$, among the~$  {\mathcal N} $ hard spheres randomly
distributed at time zero.
\end{proposition}
This generalizes  (\ref{eq: marginal time 0}) and the representation \eqref{representation}, in the sense that, for $H_n (Z_n^*([0,t] )) = h_n(Z_n^*(t))$, we obtain 
$$ F^\eps _{n, [0,t]}  ( H_n) = \int  F^\eps _n (t,Z_n^*) h_n(Z_n^*) dZ_n^*\,.$$

 \section{Dynamical correlations}
 \label{sec: Dynamical correlations}

  \subsection{Recollisions and overlaps}
  
We start from the representation \eqref{integralH} of $ F_{n, [0,t]} ^\eps(H_n )$  in terms of collision trees and pseudo-trajectories. We assume that
  $$
H_n = H^{\otimes n}
$$
with $H$ a measurable function on the Skorokhod space of trajectories $D([0,t])$ in $\D$, and we abbreviate
$$
\cH \big( \Psi_{n}^\eps  \big) :=   H^{\otimes n} \big( Z_n^*([0,t])\big) \,.
$$
Recall that there are two types of interactions between particles:
\begin{itemize}
\item a collision corresponds to the addition of a new particle;
\item recollisions occuring when two pre-existing particles collide.
\end{itemize}

\begin{figure}[h] \label{fig:rectreesf}
\centering
\includegraphics[width=4in]{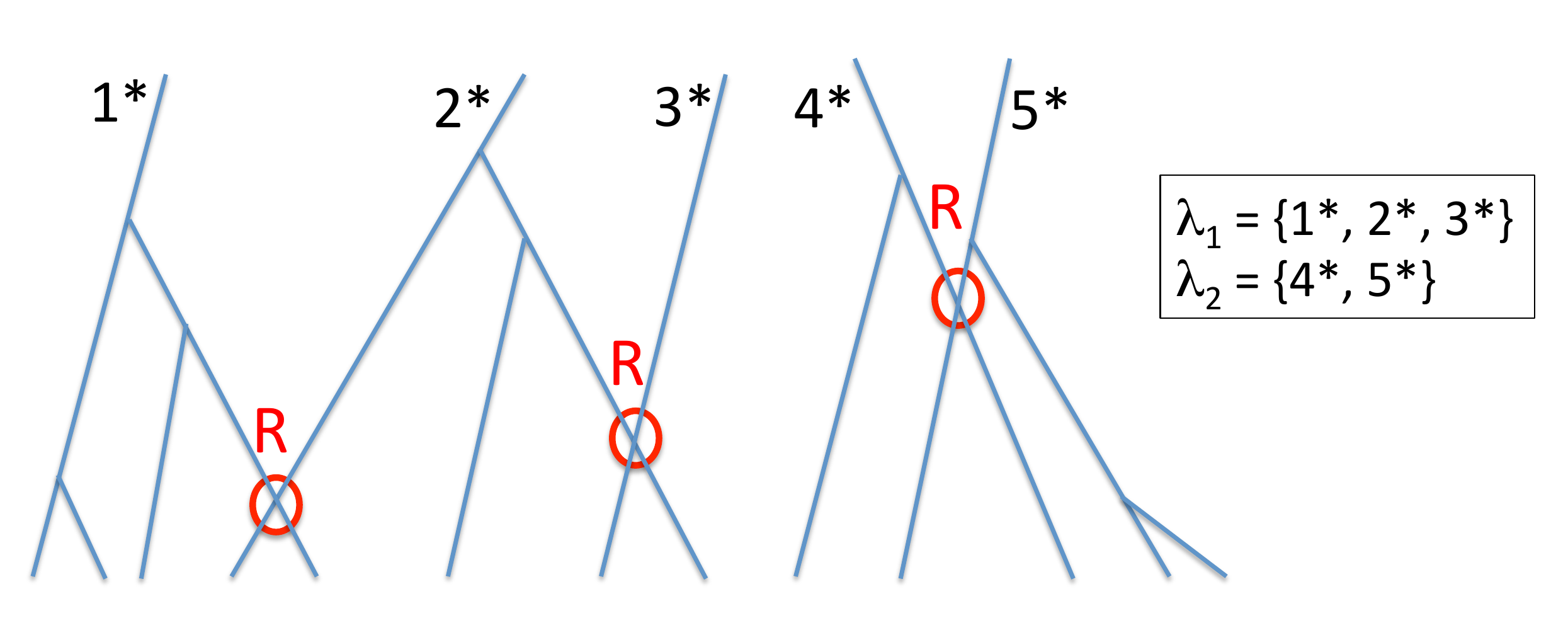} 
\caption{\small 
Recollisions connect trees into forests.
}
\label{forest}
\end{figure}

The elementary  integrals in the series expansion of $ F_{n, [0,t]} ^\eps(H_n )$ can be 
decomposed depending on whether collision trees are correlated or not by recollisions (see Figure~\ref{forest}). We then have a partition of $\{1*,\dots, n*\}$ into a certain number (say $\ell$) of {\it forests}~$(\lambda_i)_{i= 1, \dots,\ell}$, and we shall denote by~$\mb_{\gl _i} $ the characteristic function of the forest~$\lambda_i$. Namely, $\mb_{\gl _i} =1$ if and only if 
any two elements of $\lambda_i$ are connected (through their collision trees) by a chain of recollisions.
We say that $\mb_{\gl _i} =1$ is supported on {\it clusters} of size $|\lambda_i|$, formed by $|\lambda_i|$ ``recolliding'' collision trees.
We will further indicate the decomposition in forests by $\lambda = (\lambda_i)_{i= 1, \dots,\ell}$.

Formula (\ref{integralH}) can then be rewritten as a partially factorized expression:
\begin{equation} \label{eq:CE1}
\begin{aligned}
F_{n, [0,t] }^\eps(H^{\otimes n} )=  \int dZ_n^* \sum_{\ell =1}^n \sum_{\gl \in \cP_n^\ell}
 \int \cK_\lambda\left(\Psi_{\gl}^\eps\right)\,   \gP_\ell\,
 F^{\eps 0} \big (  \Psi^{\eps0}_n\big)\,    
 \end{aligned}
\end{equation}
where 
$$  \cK_\lambda\left(\Psi_{\gl}^\eps\right) =  \prod_{i=1}^{\ell } \big[ \mu (d \Psi_{\gl _i}^\eps) 
     \mb_{\gl _i}   \cC \big( \Psi_{\gl _i}^\eps \big)
  \indc_{\cG^\e} \big (  \Psi_{\gl _i}^\eps \big)\cH \big( \Psi_{\gl_i}^\eps \big) \big]\;,$$
 and $\gP_\ell= \gP_\ell \big( \gl_1, \dots, \gl_\ell \big)$ is the indicator function that particles belonging to different forests
 keep mutual distance   larger than~$\eps$.
Here and below, we indicate by~$\Psi^\eps_\a$ the pseudo-trajectory constructed starting from 
$Z_\a^*$, for any~$\a $ in~$ \{1*,\dots,n*\}$.

Although there cannot be any recollision between particles of different forests~$\lambda_i$, such particles are not yet independent, as the parameters of the pseudo-trajectories are constrained precisely by the fact that no recollision should occur. The characteristic function $\Phi_\ell = \gP_\ell \big( \gl_1, \dots, \gl_\ell \big)$ expresses this no-recollision condition. Next, we write its cumulant expansion (the analogue of \eqref{inversion}):
\begin{equation} \label{eq:CE2}
\Phi_\ell=  \sum_{r = 1}^{\ell} \sum_{\gr \in \cP^r_\ell}\varphi_\rho\,. 
\end{equation}
This formula reorganizes the $\ell$ forests into a group of $r$ {\it jungles} $\gr = \left(\gr_i\right)_{i = 1,\dots, r}$.

By construction, particles belonging to different forests will never collide among themselves. However they are allowed to ``overlap''.
We say that two different forests $\gl_i$ and $\gl_j$ {\it overlap} if two particles, belonging to the pseudo-trajectories $\Psi^\e_{\lambda_i}$ and $\Psi^\e_{\lambda_j}$ respectively, touch each other (without colliding) and cross each other freely.
Standard combinatorial arguments show then that the cumulant $\varphi_s$ of order $s$ is supported on clusters of size $s$, formed by $s$ overlapping forests (namely any two forests are connected by a chain of overlaps).

The last source of correlation in \eqref{integralH} comes from the initial data. For each given~$\rho$, we introduce a cumulant expansion of the initial data associated to $\rho$:
\begin{equation} \label{eq:CE3}
\begin{aligned}
F^{\eps 0 }\big (  \Psi^{\eps0}_n\big)  =  \sum_{s =1}^r    \sum_{\gs \in \cP_r^s}  
  f^{\eps 0}_{\gs}\,,  \quad   f^{\eps 0}_{\gs} =\prod_{i=1} ^s f^{\eps 0}_{\sigma_i} \, , \quad
  f^{\eps 0} _{\sigma_i}  =f^{\eps 0}  \left(\Psi^{\eps0}_{\sigma_i}\right)\;.
\end{aligned}
\end{equation}
Here and below, by abuse of notation, the partitions $\gs, \gr$ are also interpreted as a partition of~$\{1*,\dots, n*\} $, coarser than the partition $\gl$; the relative coarseness will be denoted by $\gl \hookrightarrow \gr \hookrightarrow \gl\, .$
 Therefore $f^{\eps 0} \big (\Psi_{\gs_i}^{\eps0}\big)$
is the time-zero cumulant evaluated on the configuration (at time 0) of 
the pseudo-trajectory starting from $Z_{\gs_i}^*$.

\begin{figure}[h] 
\centering
\includegraphics[width=4in]{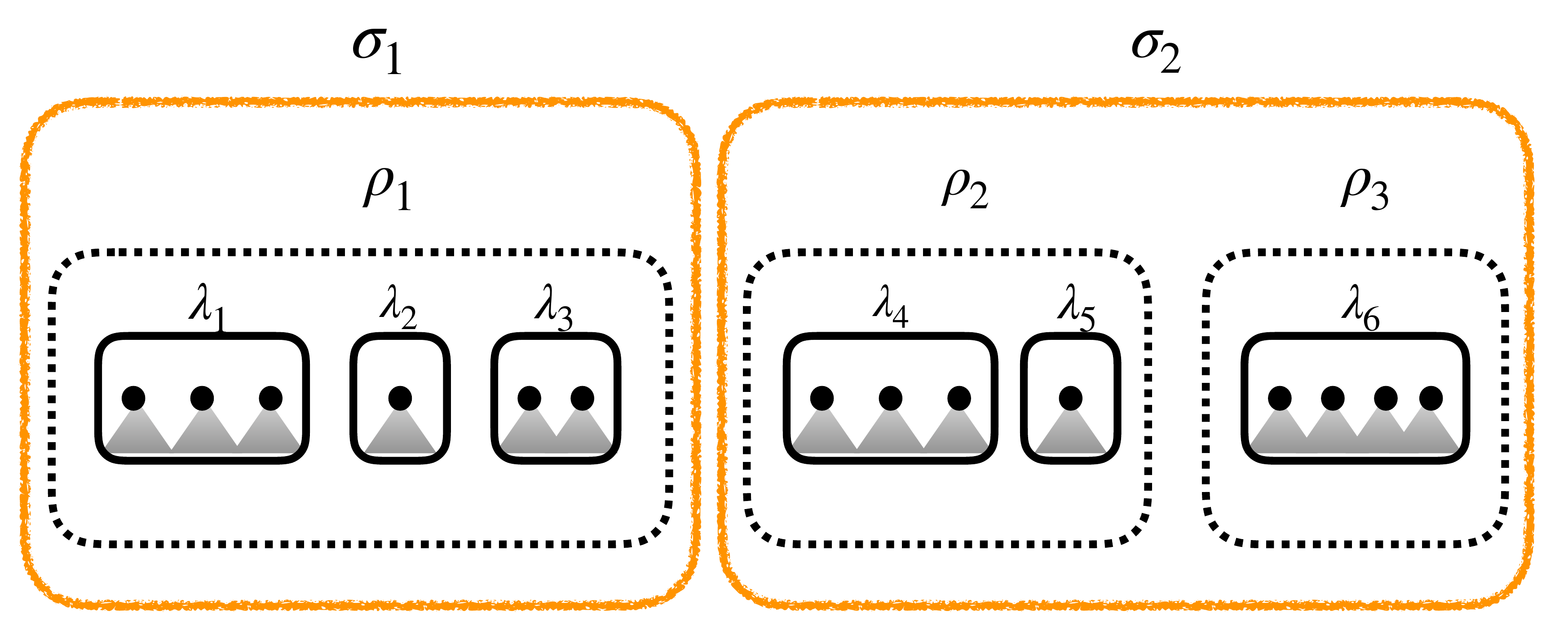} 
\caption{\small 
Clustering structure due to recollisions, overlaps and initial correlations.
}
\label{clusteringfig}
\end{figure}

%
We end up with a cluster structure on collision trees, of the form  depicted in Figure~\ref{clusteringfig}.

Replacing \eqref{eq:CE2} and \eqref{eq:CE3} into \eqref{eq:CE1}, we arrive to the following decomposition of correlation functions~:
 \begin{equation} \label{eq:CE4}
 F_{n, [0,t] }^\eps(H^{\otimes n} )  =  \int dZ_n^*\!    \sum_{\gl , \gr, \gs \atop \gl  \hookrightarrow \gr \hookrightarrow \gs}
 \int  \cK_\lambda \,
  \gp_{ \gr} \,      f^{\eps 0}_{\gs}\;,
  \end{equation}
where $\lambda$ is the partition of $\{1*, \dots , n*\}$ into $\ell$ forests of recolliding trees, $\rho$ is the partition of $\{1, \dots, \ell\}$ into $r$ jungles of overlapping forests, and $\sigma$ is the partition  of $\{1, \dots , r\}$ into initially correlated clusters. 
  
    \subsection{Cumulants and clusters}

 Comparing  formula \eqref{eq:CE4} with (\ref{inversion}), we finally identify the rescaled  dynamical cumulants (averaged over trajectories)~:
 \begin{equation}
\begin{aligned}
\label{eq: decomposition cumulant}
\cum =   \mu_\eps^{n-1} \int  dZ_n^*\sum_{\ell =1}^n \sum_{\gl \in \cP_n^\ell}
\sum_{r =1}^\ell   \sum_{\gr \in \cP_\ell^r}  
\int \cK_\lambda\,  \gp_{ \gr} \,   f^{\eps 0}_{\{1,\dots,r\}} (\Psi_{\rho_1}^{\eps 0}, \dots , \Psi_{\rho_r}^{\eps 0})  \,    .
\end{aligned}
\end{equation}
This result shows that the cumulant of order $n$ is geometrically represented by connected  clusters of size $n$~: $f_{n, [0,t] }^\eps$ corresponds  to pseudo-trajectories where the~$n$ collision trees are connected by recollisions, overlaps, or initial correlations.
  This graphical representation of cumulants leads to the following result.
  
\begin{proposition}[{\bf Convergence of dynamical cumulants}, \cite{BGSRS}] \label{prop:CDC}
Consider a  gas  of hard spheres initially distributed according to  {\rm(\ref{eq: initial measure})}.
Let  $H $  be a bounded continuous functional on $D([0,T^\star])$.
Define 
 the rescaled cumulant  $ \cum$  by {\rm(\ref{eq: decomposition cumulant})}. 
 Then, 
 \begin{itemize}
 \item
 there exists a  positive constant $C$ such that the following uniform a priori bound holds
 $$ | \cum|  \leq C^{n-1} \| H\|_\infty^n   (t+\eps) ^{n-1}   n!$$
 uniformly in $\eps$ and $n$, for any $t \leq T^\star$;
 \item when $\eps \to 0$, in the same time interval, $  \cum $ converges to a limiting $\lcum $, which is represented by a sum over minimally connected graphs, and by pseudo-trajectories with exactly $n-1$ pointwise recollisions or overlaps.
\end{itemize}
\end{proposition}

The key point  to obtain the right scaling of cumulants is to identify  ``independent''   clustering constraints~: for fixed 
$\gl, 	\gr$, collision parameters $a$ and $(T_m, \Omega_m, V_m)$ and initial velocities $V^*_n$~:
\begin{itemize}
\item
we extract a sequence of $|\lambda_i | - 1$ clustering recollisions in each forest $\lambda_i$,   a 
sequence of $|\rho_j| - 1$ clustering overlaps in each jungle $\rho_j$, and a sequence of~$r- 1$ clustering initial correlations, and prove that the factor $n!$ accounts for the combinatorics of these clustering constraints;
\item we then show that the clustering constraints can be expressed as $n-1$ conditions on the  positions at time $t$ of the particles of the pseudo-trajectory $(x_j^*)_{j=1,\dots, n}$, which are satisfied on a set of volume $O(\mu_\eps^{-(n-1)})$.
\end{itemize}
These estimates being essentially uniform with respect to the 
collision parameters~$(a, T_k, \Omega_k, V_k)$, we can sum/integrate  to 
get the  $L^1$-bound.

There is a subtle point here: a brute expansion of the overlap constraint~$\gp_s$ defined by \eqref{eq:CE2}, leads to $2^{s^2}$ terms, and cancellations need to be exploited to show that the effective number is bounded by $s!$. How to do this is known by cluster expansion techniques (see e.g.\,\cite{Ru69,GBG,PU09}). In fact, $\gp_s$ can be regarded as an Ursell function (\cite{Ru69}) by writing formally ``$\Phi_\ell(\gl_1,\cdots,\gl_\ell) = \exp\left(- U_\ell(\gl_1,\cdots,\gl_\ell)\right)$'' and interpreting $U$ as a hard core interaction on dynamical collision trees.

The proof of the second statement of Proposition \ref{prop:CDC} is very similar to Lanford's proof. We first discard the contribution of  initial correlations (which are of order $O(\eps^d)$ instead of $O(\eps^{d-1})$). 
We then prove (as discussed in Section~\ref{loops})  that any recollision which is {\it not} of clustering type will create some extra smallness, giving a vanishing contribution to the limit.

\subsection{Cumulant generating function}

The cumulants allow to characterize exponential moments of the empirical measure, as shown by the following identity~:
\begin{equation}
\label{exponential moment}
\gL^\eps_{[0,t]} (H) 
:=\frac{1}{\mu_\eps} \log \bbE_\eps \left( \exp \Big(  \sum_{i =1}^\cN H \big( {\bf z}^\eps_i ([0,t] \big)  \Big)   \right)=\sum_{n = 1}^\infty {1\over n!}  f^\eps_{n, [0,t]} \big(( e^H - 1)^{\otimes n} \big)\,,
\end{equation}
valid for functionals $H : D([0,t]) \to \bbR$ such that the series is absolutely convergent.
In order to describe the asymptotic behavior of these exponential moments, we need to obtain dynamical equations for the  limiting cumulant generating function
\begin{equation}
\label{defLambda}
 \sum_{n = 1}^\infty {1\over n!}  f_{n, [0,t]}\big(( e^H - 1)^{\otimes n} \big) \,,
\end{equation}
which is well defined (as a corollary of Proposition \ref{prop:CDC}) for $t\in [0,T^\star]$ provided that $H$ is a continuous functional
satisfying a suitable bound:
$$  \Big| \Big (e^{H \big (z([0,t]) \big)}-1\Big)^{\otimes n}   \Big|  \leq \exp\Big(   \alpha_0  n +\frac{\beta_0}4\sup_{s\in [0,t] }   |V_n(s)|^2\Big)\,,$$
for some $\alpha_0$ (related to the constant $C_0$ in \eqref{beta0} and to~$T^\star$).

We shall not write here the hierarchical equations for the family of cumulants at equal times $\left(f_{n}(t)\right)_{n \geq 1}$, obtained by choosing $H (z([0,t] )) = h(z(t))$ in~\eqref{exponential moment}. This hierarchy (mentioned in Section \ref{strategy} as  ``Boltzmann cumulant hierarchy'') is derived and analysed in \cite{EC81}. Our purpose is to focus directly on the full series \eqref{defLambda}, which we study for a class of regular test functionals.

For $t\in [0, T^\star]$, denote by~$ \cJ(t, \varphi, \gamma)$ the limiting cumulant generating function~(\ref{defLambda}) associated with
 $$
e^{ H \big (z([0,t]) \big)} = \gamma  \big (z(t)\big) \exp \Big(  -  \int_0^t  \gp  \big ( s,z(s) \big) ds 
\Big)\, ,
$$
where   $(\gp ,\gamma )$  belong to 
\begin{equation*}
\begin{aligned}
{\mathcal B} :=    &\Big\{(\varphi,\gamma) \in  C^1([0,t] \times \D;{\mathbb C}) \times C^1(\D;{\mathbb C} )
\, \   \big|\  \, \\  & | \gamma(z)| \leq  e^{\frac{1}{2}( \alpha_0 +\frac{\beta_0}{4} |v|^2)} ,\quad
    \sup_{s\in [0,T^\star]}| \gp (s,z) | \leq {1\over 2T^\star} \left(\alpha_0  +\frac{\beta_0}{4} |v|^2\right)
\Big\}\, .
\end{aligned}
\end{equation*}
We shall be interested in functions of the form
 $$
 \varphi  = D_sh \equiv (\d_s+ v\cdot \nabla_x) h\quad \hbox{ and } \quad \gamma  = \exp (h(t)) \,,$$
 therefore we simplify notation by setting~$$ \cJ (t,h) :=  \cJ(t, Dh, \gamma)_{| \gamma = \exp (h(t))} \;.$$
We set $\bbB := \left\{\, h  \,\   \big|\  \,  (D_t h, \exp( h(T^\star)) ) \in {\mathcal B}\,\right\}\,.$
With these notations, the following result holds.

\begin{theorem}[{\bf Hamilton-Jacobi equations}, \cite{BGSRS}]
\label{thm: HJ}
The  functional $\cJ$ is analytic with respect to $\g$  on $\bbB$, and it satisfies on $[0,T^\star]$ the Hamilton-Jacobi equation 
\begin{align}
\label{eq: HJ}
\d_t \cJ (t,h)  & = \frac12  \int {\d \cJ (t,h) \over \d \gamma} (z_1)  {\d \cJ(t,h) \over \d \gamma } (z_2) \\
&\qquad \times \Big (e^{h(t,z'_1) +h(t,z'_2)} - e^{h(t,z_1) +h(t,z_2)}\Big )   d\mu (z_1, z_2, \omega)\, , 
\nonumber
\end{align}
where 
\begin{equation}\label{defdmu}
d\mu(z_1, z_2, \omega):=  \delta({x_1 - x_2})\,  ((v_1-v_2)\cdot \omega)_+ d\omega dv_1 dv_2 dx_1 dx_2\,.\end{equation}
\end{theorem}

%
%

The local existence and uniqueness of a solution for this Hamilton-Jacobi equation relies on a Cauchy-Kowalewski argument in a functional space,  encoding the loss continuity estimates due to the divergence of the collision cross section~(\ref{defdmu}) at  large velocities.

\section{The fluctuating Boltzmann equation}

Describing the fluctuations around the Boltzmann equation  is a first way to 
capture part of the  information which has been lost in the limit~\eqref{eq:R1}. 
As in the standard central limit theorem, we expect these fluctuations to be of order~$1/\sqrt{\mu_\eps}$. 
We therefore define  the fluctuation field $\gz^\gep$ by (see \ref{eq: fluctuation field})
$$
\gz^\gep_t \big(  h  \big) :=  { \sqrt{\mu_\eps }} 
\left( \pi^\eps_t(h) -  \int \, F^\eps_1(t,z) \,  h \big(  z \big)\, dz   \right) \, ,
$$
 for any test function $h: \D\to\R$.
 
It is easy to check that, in our assumptions, the empirical measure starts close to the density profile $f^0$ and that~$\gz_0^\eps$ converges to a Gaussian white noise~$\zeta_0$ with covariance $$\bbE \left( \gz_0(h_1)\, \gz_0(h_2) \right)
 = \int h_1(z) \,h_2(z)\, f^0(z) \,dz\, .$$

Moreover, it follows from Proposition \ref{prop:CDC} that
$\zeta^\eps$ converges  to a solution of the  fluctuating Boltzmann equation 
\begin{equation}
\label{eq: OU}
d \gz_t   = \cL_t \,\gz_t\, dt + d\eta_t\,,
\end{equation}
where $\cL_t $ is the  {linearized Boltzmann operator} around the solution $f$ of the Boltzmann equation~\eqref{boltzmann-eq} 
$$\begin{aligned}
&\cL_t \,h(x,v) := - v \cdot \nabla_x h(x,v)+\int_{\R^d}\int_{{\mathbb S}^{d-1}} \,d \omega \, d w   \left( (v - w) \cdot \omega \right)_+ \\
&     \times \big( f (t,x,w') h(x,v') + f (t,x,v') h(x,w') 
   - f (t,x,v)  h(x,w) -  f (t,x,w)  h(x,v) \big),
\end{aligned}
$$
and $d \eta_t(x,v)$ is a Gaussian noise  with zero mean and covariance 
\begin{align}
\label{eq: covariance bruit}
& \bbE \left( \int  dt_1\, dz_1 \,  h_1 (z_1)  \,\eta_{t_1} (z_1) \, \int dt_2 \, dz_2 \, h_2 (z_2) \,\eta_{t_2}(z_2) \right)
  \\
  &\qquad\qquad = \frac{1}{2} \int  dt  \int 
 f (t, z_1)\, f (t, z_2) \Delta h_1 \, \Delta h_2 \; d\mu(z_1, z_2, \omega) \nonumber
\end{align}
with notation \eqref{defdmu} and
\begin{equation}\label{defDelta}
\Delta h (z_1, z_2, \omega) := h(x_1,v'_1) +  h(x_2,v'_2) -  h(x_1,v_1) -  h(x_2,v_2) \,.
\end{equation}

Our main result is   then the  following. 

\begin{theorem} [{\bf Fluctuating Boltzmann equation}, \cite{BGSRS}]
\label{thmTCL}
Consider a system of hard spheres initially distributed according to  {\rm(\ref{eq: initial measure})}. 
Then, in the Boltzmann-Grad limit $\mu_\eps \to \infty$, the fluctuation field $\left(\zeta^\eps_t\right)_{t \geq 0}$ converges in law on $[0,T^\star]$ to the solution~$\left(\zeta_t\right)_{t \geq 0}$ of the   fluctuating Boltzmann equation~{\rm(\ref{eq: OU})}.
\end{theorem}

The convergence towards the limiting process \eqref{eq: OU} was conjectured  by Spohn in \cite{S83} and 
the non-equilibrium covariance of the process at two different times was obtained in \cite{S81}. 
The noise emerges after averaging the deterministic microscopic dynamics.
It is white in time and space, but correlated in   
velocities so that momentum and energy are conserved.

We further recall a few properties (referring to \cite{EC81,S81,S83} for details).

\begin{itemize}
\item In the equilibrium case ($f^0(x,v)= M_{\beta}(v)$ where $M_{\beta}$ is a Maxwellian with inverse temperature $\b$) the noise term   compensates the dissipation induced by  the (stationary) linearized collision operator
$\cL$, and the covariance of the noise  can be predicted heuristically by using the invariant measure.

\item 
Out of equilibrium, on the one hand, the noise covariance \eqref{eq: covariance bruit} can be simply understood as a generalization of the covariance at equilibrium, based on the assumption (which can be proved for short times \cite{S2}) that the system is locally Poisson distributed; i.e. on a small cube around $x$ at time $t$ we see a uniform ideal gas with density $\int f(t, x,v) dv$ and velocity distribution $f (t,x,v) / \int f(t, x,v_*) dv_*$. The noise being delta-correlated in space and time, its structure is obtained from the equilibrium case after the replacement $M_\b(v) \to f(t,x,v)$.
 
 \item On the other hand, 
the covariance of the fluctuation field out of equilibrium has a subtle microscopic structure 
originating from recollisions in the Newtonian dynamics.
To see this, it is enough to compute the covariance of the fluctuation field at time $t$ by using \eqref{Eepsh2}~:
\begin{align*}
 \bbE_\eps \Big(  \gz^\gep_t \big(  h  \big)^2 \Big)
 = &   \int F^\eps_1(t,z_1)\, h^2(z_1)\, d z_1  \\
& \qquad +  \int \mu_\eps \Big( F^\eps_2(t,Z_2) - F^\eps_1(t,z_1) F^\eps_1(t,z_2) \Big)\, h(z_1) h(z_2)\, d Z_2 \\
= & \int f^\eps_1(t,z_1)\, h^2(z_1)\, d z_1 +  \int f^\eps_2(t,Z_2) \, h(z_1) h(z_2)\, d Z_2 
\end{align*}
where, in the second equality, we used the first two cumulants as defined by \eqref{inversion}.
The last term is zero at equilibrium, while out of equilibrium describes correlations visible at macroscopic distance in space.
But, as made apparent from the geometrical representation \eqref{eq: decomposition cumulant}, $f^\eps_2 $ records the effect of one (and only one) recollision/overlap; meaning precisely that the pseudo-trajectories contributing to $f^\eps_2 $ have the form in Figure \ref{noise-fig}.
 \begin{figure}[h] 
\centering
\includegraphics[width=4in]{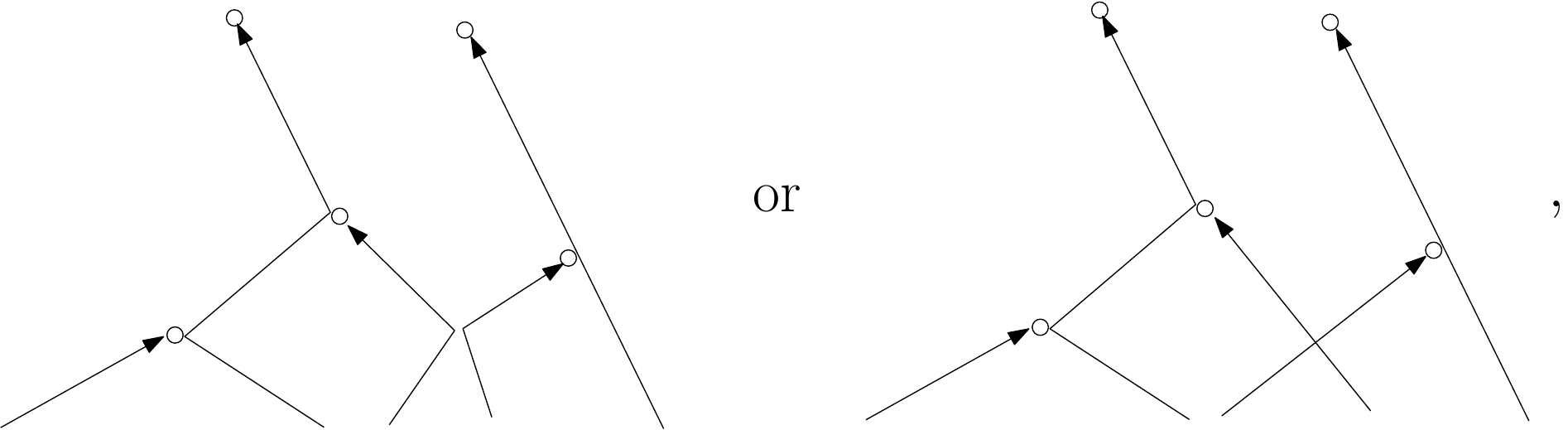} 
\caption{Trajectories contributing to the equal time covariance at two different points in space (case of $m=3$ collisions with fresh particles): clustering recollisions or clustering overlaps}
\label{noise-fig}
\end{figure}
Contrary to the typical behavior of the hard sphere gas for which recollisions can be  neglected, 
the  covariance of the limiting Gaussian process encodes exactly the effect of a single recollision.
\end{itemize}

The uniform bounds on the cumulants discussed in the previous section are considerably better than what is  required to obtain Theorem \ref{thmTCL}. The proof amounts indeed to looking at a characteristic function living on larger scales. 
A more technical part concerns the tightness of the process. This can be achieved adapting a Garsia-Rodemich-Rumsey's inequality on the modulus of continuity, to the case of a discontinuous process. We omit the details, and focus on the characteristic function only.

  Consider the function $H$ defined by 
\begin{equation}
\label{discrete-H}
H (z([0,t])) = \sum_{p=1} ^P  h_p \big(z (\theta_p)\big)  
\end{equation}
for a  finite sequence of times~$(\theta_p )_{1\leq p \leq P} $ and  weights~$(h_p )_{1\leq p \leq P} $.
The characteristic function can be rewritten  in terms of the empirical measure
$$
\begin{aligned}
\log \bbE_\eps  \left( \exp  \Big(  \sum_{p= 1} ^P \zeta^\eps_{\theta_p} (h_p)  \Big) \right)   
&= \mu_\eps \sum_{ n=1}^\infty   \frac{1}{n !}
 f_{n, [0,t] }^\eps \left( \big( e^{  {H\over \sqrt{\mu_\eps}}} - 1\big)^{\otimes n} \right)  \\
 &\qquad 
-  \sqrt{\mu_\eps}     \sum_{p= 1}^P  \int  F^\eps_1(\theta_p,z)  h_p(z)\,dz \, .
\end{aligned}
$$
At leading order, only the terms $n=1$ and $n=2$ will be relevant  in the limit since
$$
\Big| f_{n, [0,t] }^\eps \left( \big( e^{ H\over \sqrt{\mu_\eps}} - 1\big)^{\otimes n} \right)  \Big| 
\leq 
C^n \left\| e^{ H\over \sqrt{\mu_\eps}} - 1\right\|_\infty ^n n! \, .
$$
Expanding the exponential  with respect to $\mu_\eps$, we also notice that the  term 
of  order~$\sqrt{\mu_\eps}$ cancels and we find
$$
\log \bbE_\eps  \left( \exp  \Big( \sum_{p= 1} ^P \zeta^\eps_{\theta_p} (h_p)  \Big) \right)    
 =  \frac{1}{2}     f_{1, [0,t]} ^\eps  \left( H^2 \right)
 + \frac{1}{2}  f_{2, [0,t]}^\eps  \left( H^{\otimes 2} \right)
+ O \left( \frac{1}{\sqrt{ \mu_\eps} } \right) \,.$$
  Then  the characteristic function $\bbE_\eps  \left(  \exp \big(  \sum_{p= 1} ^P \zeta^\eps_{\theta_p} (h_p)   \big) \right)$
converges to the characteristic function of a Gaussian process. 

From the equations on $f_1$ and $ f_2$, we deduce that  the limiting covariance $\cC = \cC(s,t,\gp,\psi)$ satisfies the following dynamical equations, for test functions $\gp,\psi$ on $ \D$~:
$$\begin{cases}
\d_t  \cC(s,t,\gp,\psi)    =  \cC(s,t,\gp, \cL_t^*\psi)\,, \\
\d_t  \cC(t,t,\gp,\psi)    =  \cC(t,t,\gp, \cL_t^*\psi) +  \cC(t,t, \cL_t^* \gp, \psi)
+ {\bf  Cov}_t (    \psi   ,  \gp  )\,, \\
\cC(0,0,\gp,\psi) =\displaystyle  \int \gp(z) \psi(z) f^0(z) dz\,,
\end{cases}
$$
where   
$$
{\bf  Cov}_t ( \gp, \psi) 
:= \frac{1}{2}  \int d\mu (z_1, z_2, \omega) \,
 f (t, z_1)\, f (t, z_2)   \Delta \psi \Delta \gp$$
 with notation~(\ref{defdmu}) and~(\ref{defDelta}),
and $ \cL_t^* := v\cdot \nabla_x + {\bf L}_t^*$  with 
$$
{\bf L}_t^* \,\gp(v) := \int_{\R^d}\int_{{\mathbb S}^{d-1}}  \,d \omega \, d w   \left( (v - w) \cdot \nu \right)_+ 
 f (t,w)  \Delta \gp \, .
$$
The covariance of the  fluctuating Boltzmann equation (\ref{eq: OU}) satisfies the same equations, and we conclude by a uniqueness argument  that both processes coincide.

 \section{Large deviations}
 
 While typical fluctuations are of order $O(\mu_\eps^{-1/2})$, 
 larger fluctuations may sometimes happen, leading to  an evolution which is different from  the typical one given by  the Boltzmann equation.
A classical problem is  to evaluate the probability of such atypical events, namely that the empirical measure $\pi^\eps_t$, defined in~\eqref{eq: empirical},
remains close to a probability density $\gp_t$   during the time interval~$[0,T^\star]$.

In the G\"artner-Ellis theory of large deviations \cite{dembozeitouni}, 
the large deviation functional is given as the Legendre transform of the limiting cumulant 
generating function. 
The outcome of the cumulant analysis was the existence of the 
limiting exponential moment~$\cJ(t,h)$ and its characterization 
via the Hamilton-Jacobi equation in Theorem~\ref{thm: HJ}.
For any $t \leq T^\star$, we then define the large deviation functional on the time interval $[0,t]$ as  
\begin{equation}
\label{eq: LD cF}
\begin{aligned}
\cF (t, \varphi) &:= \sup_{ h \in \bbB} \Big\{   - \int_0^t \int_{\mathbb D} \varphi(s, z) D_s h (s, z) dz ds  \\
&\qquad\qquad + \int_{\mathbb D} \varphi(t,z) h(t,z) dz  
-  \cJ(t,h)  \Big\} \, .
\end{aligned}
\end{equation} 
Since the supremum is restricted, for technical reasons, to the test functions in~$\bbB$, we do not expect   
$\cF$ to be the correct large deviation functional.
However the following theorem shows that the functional $\cF$ fully describes the 
large deviation behavior for densities $\gp$ such that the supremum in \eqref{eq: LD cF}
is reached for some $h \in \bbB$. This restricted set of densities $\gp$ will be called $\cR$.

A different, explicit formula for the large deviation functional was obtained by Rezakhanlou \cite{Rez2} in the case of a one-dimensional stochastic dynamics mimicking the hard-sphere dynamics, and then 
conjectured for the three-dimensional, deterministic hard-sphere dynamics by Bouchet \cite{bouchet}~:
\begin{align}
&\widehat\cF(t,\gp) 
:=  \widehat \cF (0,\gp_0)  \nonumber \\
& + 
\sup_p \left\{ \int_0^{t} ds \left[
\int_{\bbT^d} dx  \, \int_{\bbR^d} dv \,  p(s,x,v)  \, D_s\gp(s,x,v) - \cH \big( \gp(s) ,p(s) \big) \right]
\right\} ,
\label{eq: hat cF}
\end{align}
where the supremum is taken over bounded measurable functions $p$ growing at most quadratically in $v$,  the Hamiltonian is given by
$$\cH(\gp,p) := \frac{1}{2} \int  d\mu (z_1, z_2, \omega) 
\gp (z_1) \gp (z_2) \big( \exp \big( \Delta p \big) -1 \big) 
$$
and $\widehat \cF(0,\cdot)$ stands for the large deviation functional on the initial data
\begin{equation}
\widehat \cF (0,\gp_0) 
= \int dz   \left( \gp_0 \log \left( \frac{\gp_0}{f^0} \right) -  \gp_0 + f^0 \right).
\end{equation}
Let $\hat \cR$ denote the set  of densities $\gp$ such that the supremum in \eqref{eq: hat cF}
is reached for some $p \in \bbB$.

Let $\cM(\bbD)$ be the set of probability measures on $\bbD$.

Our main result is then the following.

\begin{theorem} [{\bf Large deviations, \cite{BGSRS}}]
\label{thmLD}
Consider a system of hard spheres initially distributed according to   {\rm(\ref{eq: initial measure})}. 
In the Boltzmann-Grad limit $\mu_\eps \to \infty$, the empirical measure $\pi^\eps$ satisfies the following large deviation estimates  for any $t \in [0,T^\star]$.
\begin{itemize}
\item For any compact set ${\bf F}$ of the Skorokhod space $D([0,T^\star], \cM)$, 
\begin{align}
\label{eq: large deviations upper bound}
\limsup_{\mu_\eps \to \infty} \frac{1}{\mu_\eps} 
\log \bbP_\eps \left( \pi^\gep  \in {\bf F} \right) \leq - \inf_{\gp \in \bf F}  \cF(T^\star,\gp)\, .
\end{align}
\item For any open set ${\bf O}$ of the Skorokhod space $D([0,T^\star], \cM)$, 
\begin{align}
\label{eq: large deviations lower bound}
\liminf_{\mu_\eps \to \infty} \frac{1}{\mu_\eps} 
\log \bbP_\eps \left( \pi^\gep  \in {\bf O} \right) \geq - \inf_{\gp \in \bf O \cap \cR} \cF(T^\star,\gp)\, .
\end{align}
\end{itemize}
Moreover, for any $\gp \in \cR \cap \hat \cR$ and $t$ sufficiently small, one has that
$\cF(t,\gp) = \widehat \cF(t,\gp)$.
\end{theorem}

Given our precise control of the exponential moments, the large deviation proof is standard.
Note that, in absence of global convexity, we cannot succeed in proving a full large deviation principle. However, restricting to a class of regular profiles, the variational problem defining the dual of $\widehat \cF$
 can be uniquely solved and identified with the solution of the Hamilton-Jacobi equation \eqref{eq: HJ}. 
 The result then follows from a uniqueness property of \eqref{eq: HJ}.

 \bigskip

\noindent {\bf Acknowledgements} We are very grateful to H. Spohn, F. Bouchet, F. Rezakhanlou,   G. Basile, D. Benedetto, and L. Bertini for many enlightening discussions on the subjects treated in this text.
T.B. acknowledges the support of ANR-15-CE40-0020-01 grant LSD.

\end{document}